\documentclass{article}%
\usepackage{amssymb}
\usepackage{amsmath}
\usepackage{amsfonts}
\usepackage{graphicx}%
\providecommand{\U}[1]{\protect\rule{.1in}{.1in}}
\newtheorem{theorem}{Theorem}

\newtheorem{lemma}[theorem]{Lemma}

\newtheorem{proposition}[theorem]{Proposition}

\begin{document}

\title{\ A geometric characterization of Poisson type distributions}
\author{V. P. Palamodov}
\date{}
\maketitle

\textbf{Abstract: }Tempered distributions on $\mathbb{R}^{n}$ that have the
Poisson\ structure are characterized in terms of geometry of its support and spectrum



\section{Introduction}

It is known in crystallography that the diffracted wave function is the
Fourier transform of the electronic density function. The old crystallographic
wisdom was whenever one sees isolated Bragg's peaks in the Fourier spectrum,
one must have a periodic structure, i.e. a crystal. The term \textquotedblleft
quasicrystal\textquotedblright\ appeared after the discovery of "a new
principle of packing of molecules" \cite{S}. This term is used in mathematics
as the name of distributions that has closed discrete support and diffraction
spectrum (support of the Fourier transform) and some additional properties. A
survey of the theory of quasicrystals is given in Lagarias \cite{La00}. A set
$\Lambda\subset\mathbb{R}^{n}$ is called \textit{uniformly discrete} (u.d.)
if
\[
d\left(  \Lambda\right)  =\mathrm{\inf}\left\vert p-p^{\prime}\right\vert
>0,\ p\neq p^{\prime}\in\Lambda.
\]
The following important result was stated by N. Lev and A. Olevskii
\cite{LO15}, conjectured by Lagarias \cite{La00}.

\begin{theorem}
\label{LO} \textit{If} $\mu=\sum_{\Lambda}\mu\left(  p\right)  \delta
_{p},\ \mu\left(  p\right)  >0$\textit{\ is a positive measure on }%
$\mathbb{R}^{n}$\textit{\ supported by a u.d.\ set }$\Lambda$ \textit{such
that the Fourier transform is also a measure }$\hat{\mu}=\sum_{\Sigma}\hat
{\mu}\left(  \sigma\right)  \delta_{\sigma},$ $\hat{\mu}\left(  \sigma\right)
\neq0$ \textit{supported on a u.d. set }$\Sigma$, \textit{then }$\Lambda
$\textit{\ is contained in the finite union of shifts of a }full rank
\textit{lattice }$L$ \textit{in }$\mathbb{R}^{n}$\textit{. The same is true
for }$\Sigma$\textit{\ and the dual lattice }$L^{\ast}$\textit{. In the case
}$n=1$ \textit{these conclusions hold without the positivity assumption.}
\end{theorem}

Previous results in this direction were obtained by C\'{o}rdoba \cite{Cor},
Kolountzakis-Lagarias \cite{Kol}, Meyer \cite{Me70} and Favorov \cite{Fav}. S.
Favorov has shown that a u.d. set $\Lambda$ must be the union of shifts of a
finite number of (may be not parallel) lattices provided that $\ \Sigma$ is
countable and the set $\left\{  \left\vert \mu_{p}\right\vert ,\ p\in
\Lambda\right\}  $ is finite. These results are based on Helson-Cohen's theory
of idempotent measures on locally compact abelian groups. This theory depends
on the choice axiom (or Tychonoff's theorem). Constructive proofs of these
results are desirable. The arguments of \cite{LO15}\ are constructive and
incorporate the concept of Meyer sets \cite{Me70}, \cite{Me95}.

A full rank lattice in $E=\mathbb{R}^{n}$ is the set $L=g\mathbb{Z}^{n},$
where $g$ is a nonsingular $n\times n$ matrix.$\ $It is the additive subgroup
of $E$\ generated by columns of $g$ (called generators of the lattice). The
lattice $L^{\ast}=h\mathbb{Z}^{n},$ $h=\left(  g^{t}\right)  ^{-1}$ in the
dual space $E^{\ast}$ is called dual. A $L$-crystal is the set $\Lambda
=\cup_{1}^{K}L+q_{k},$ where $q_{1},...,q_{K}\in E.$ S. Favorov \cite{Fav}
gave an example of a measure $\mu$ on $\mathbb{R}^{2}$ whose support $\Lambda$
is the u.d. union of\ two non parallel crystals $\Lambda_{1}$ and $\Lambda
_{2}$ and the same is true for the spectrum. Favorov's measure $\mu$ has no
Poisson structure. The purpose of this paper is to prove

\begin{theorem}
\label{Gen}If a tempered distribution $t\neq0$ on $\mathbb{R}^{n}\ $is\textit{
}supported by a set $\Lambda$ such that $\Lambda-\Lambda$ is u.d. and the
spectrum of $t$ is supported by a u.d. set $\Sigma,$ then $\Lambda$ \textit{is
the }$L$-\textit{crystal with a lattice }$L$ \textit{and }$\Sigma$ is the
$L^{\ast}$-crystal with the dual lattice $L^{\ast}$.\textit{\ }
\end{theorem}

Note that any crystal $\Gamma$ is a u.d. set and is true for the set
$\Gamma-\Gamma$. Therefore the condition $\Lambda-\Lambda$ and $\Sigma$ are
u.d. in Theorem \ref{Gen} is necessary for the conclusion. In view of symmetry
of support and spectrum of a distribution, the condition $\Lambda$ and
$\Sigma-\Sigma$ are u.d. is also a characterization of crystals. We shall also
show that any distribution $t$ as above has Poisson structure (\ref{16})
(Poisson's comb).

For measures $t=\mu$ these results were stated by Lev and Olevskii
\cite{LO13}, see in \cite{LO15} for father results. Our proof is based on the
result of \cite{LO15}. The assumption that $\Lambda-\Lambda$ is u.d. implies
that $\Lambda$ is itself a u.d. set. This later condition is not sufficient,
which follows from Favorov's example \cite{Fav}. See \cite{Mo} for a survey of
properties of sets $\Lambda$ such that $\Lambda-\Lambda$ is u.d.

I am grateful to A. Olevskii for attracting my attention to this problem and
bibliographic commentaries. I thank S. Favorov for commentaries to his papers.

\section{Preliminaries}

Let $x=\left(  x_{1},...,x_{n}\right)  $ be a system of linear coordinates in
$E\cong\mathbb{R}^{n}.$ A Schwartz (test) function on $E$ is a smooth function
$\varphi$ such that $P\left(  x\right)  D^{i}\varphi\left(  x\right)  $ is
bounded for any polynomial $P$ and any partial derivative $D^{i}%
=\partial_{x_{1}}^{i_{1}}...\partial_{x_{n}}^{i_{n}},$ $i=\left(
i_{1},...,i_{n}\right)  \in\mathbb{Z}^{n}.$ The vector space $\mathrm{S}%
\left(  E\right)  $ of all Schwartz functions is called the Schwartz space and
is supplied with the natural topology. A tempered distribution $t$ on $E$ is a
linear continuous functional on the Schwartz space see more detail in
\cite{Tr}. For a point $p\in E,$ the functional$\ \delta_{p}\left(
\varphi\right)  =\varphi\left(  p\right)  \ $is called delta distribution at
$p.$ Any partial derivative of $\delta_{p}$ is a tempered distribution. A
tempered function is a linear continuous functional on the Schwartz space of
densities $\varphi\mathrm{d}x,$ where $\varphi$ is a test function and
$\mathrm{d}x=\mathrm{d}x_{1}\wedge...\wedge\mathrm{d}x_{n}.$ For a tempered
function $f,$ the product $f\mathrm{d}x$ is a tempered distribution. The
Fourier transform $\hat{t}=F_{x\rightarrow\xi}\left(  t\right)  $ of a
tempered distribution $t$ is a tempered function on the dual space $E^{\ast}$
(L. Schwartz's theorem, see \cite{Tr}). It is defined by$\ $%
\[
\hat{t}\left(  \rho\right)  =t\left(  F\left(  \rho\right)  \right)
,\ F\left(  \rho\right)  =\int\exp\left(  -\mathrm{j}\left\langle
\xi,x\right\rangle \right)  \psi\left(  \xi\right)  \mathrm{d}\xi,
\]
where $\rho=\psi\mathrm{d}\xi$, $\psi$ is a test function on $E^{\ast}$,
$\ \xi_{1},..,\xi_{n}$ are\ dual coordinates in $E^{\ast}$ and $\mathrm{j}%
=2\pi i.$ Let $B\left(  \rho\right)  $ be the open ball in $E$ or in $E^{\ast
}$ of radius $\rho$ with the center at the origin$.$

\begin{proposition}
\label{TD}Any tempered distribution $t$ on $E$ supported on a u.d. set
$\Lambda$ is uniquely represented in the form%
\begin{equation}
t=\sum_{p\in\Lambda}t_{p}\left(  D\right)  \delta_{p}, \label{14}%
\end{equation}
where $t_{p}\left(  D\right)  $ is a differential operator on $E$ and for any
test function $\alpha$ on $E,$
\[
t\left(  \alpha\right)  =\sum_{p\in\Lambda}\sum_{\left\vert i\right\vert \leq
m\left(  p\right)  }\delta_{p}\left[  t_{p}\left(  -D\right)  \alpha\left(
x\right)  \right]  .
\]
The sum (\ref{14}) is a tempered distribution if and only if $\mathrm{ord\ }%
t\doteqdot\mathrm{sup}$\ $\mathrm{ord\ }t_{p}<\infty,$ and
\begin{equation}
\sum_{p\in\Lambda}\left(  \left\vert p\right\vert +1\right)  ^{-m}\left\Vert
t_{p}\left(  D\right)  \right\Vert <\infty\label{4}%
\end{equation}
for a constant $k,$ where
\[
\left\Vert t_{p}\left(  D\right)  \right\Vert \doteqdot\max_{i}\left\vert
t_{p,i}\right\vert ,\ t_{p}\left(  D\right)  =\sum t_{p,i}D^{i},\ t_{p,i}%
\in\mathbb{C}.
\]

\end{proposition}

\textit{Proof}. By \cite{Tr}, $t$ satisfies%
\begin{equation}
\left\vert t\left(  \varphi\right)  \right\vert \leq C\mathrm{sup}_{\left\vert
i\right\vert \leq k}\mathrm{sup}_{x\in E}\left(  1+\left\vert x\right\vert
\right)  ^{m}\left\vert D^{i}\varphi\left(  x\right)  \right\vert \label{11}%
\end{equation}
for some $k,m$ and $C>0.$ For any point $p\in\Lambda,\ $we choose a
neighborhood $U_{p}$ of $p$ such that $U_{p}\cap U_{q}=\varnothing$ for $p\neq
q.$ Let $e_{p}$ be a test function on $E$ whose support is contained in
$U_{p}$ that is equal to $1$ on a neighborhood of $p$. We have%
\[
t\left(  \varphi\right)  =\sum_{p\in\Lambda}e_{p}t\left(  \varphi\right)  ,
\]
since $\varphi-\sum e_{p}\varphi$ vanishes on a neighborhood of $\Lambda.$ The
distribution $e_{p}t$ is supported by $p.$ By L. Schwartz's theorem,
$e_{p}t=t_{p}\left(  D\right)  \delta_{p}$ for a differential operator
$t_{p}\left(  D\right)  $ see \cite{Tr}, Theorem 24.6. Inequality (\ref{11})
implies that $\mathrm{ord\ }t\leq k$ and (\ref{4}). $\square$

\section{Starting the proof}

According to \cite{Wie} and \cite{Hof1}, the convolution $t\ast t^{\ast}$ is
called the \textit{autocorrelation} of the "signal" $t,$ where
\[
t^{\ast}\left(  x\right)  =\bar{t}\left(  -x\right)  =\sum t_{i}\left(
-x\right)  \left(  -D\right)  ^{i}.
\]
Following \cite{LO15}, we define the regularized autocorrelation $A_{\varphi
}\doteq\varphi t\ast t^{\ast}$ for an arbitrary tempered $t$ taking for
$\varphi$ an arbitrary test function on $E$.

\begin{lemma}
\label{C} $\Sigma$ is a $\Gamma$-crystal with a lattice $\Gamma.$
\end{lemma}

\textit{Proof of Lemma. }We have
\[
A_{\varphi}=\sum_{p\in\Lambda}\varphi\left(  p\right)  t_{p}\left(  D\right)
\delta_{p}\ast\sum_{q\in\Lambda}\bar{t}_{q}\left(  -D\right)  \delta_{-q}.
\]
We integrate by parts near each point $q\in\Lambda$ and apply equation
$\delta_{p}\ast\delta_{-q}=\delta_{p-q}$, which yields%
\[
A_{\varphi}=\sum_{p}\sum_{q}\bar{t}_{q}\left(  -D\right)  \varphi\left(
q\right)  t_{p}\left(  D\right)  \delta_{p-q}\left(  x\right)  =\sum
_{r\in\Lambda-\Lambda}a_{\varphi,r}\delta_{r}\left(  x\right)  ,
\]
where the series
\[
a\doteqdot\sum_{q,q+r\in\Lambda}\bar{t}_{q}\left(  -D\right)  \varphi\left(
q\right)  t_{q+r}\left(  D\right)
\]
converges for any $r\in\Lambda-\Lambda$, since by Proposition \ref{TD}, the
norm $\left\Vert t_{q}\overline{t}_{q+r}\right\Vert $ has at most polynomial
growth as $q\rightarrow\infty$. This\ implies that distribution%
\[
A_{\varphi}=\sum_{r\in\Lambda-\Lambda}a_{\varphi,r}\delta_{r}\left(  x\right)
\]
is supported by the u.d. set $\Lambda-\Lambda.$ Farther, we have$\ \overline
{\hat{T}}$
\[
F_{x\mapsto\xi}\left(  A_{\varphi}\right)  =\widehat{\varphi t}\cdot
\overline{\hat{t}}=\left(  \hat{\varphi}\ast\hat{t}\right)  \cdot
\overline{\hat{t}}.
\]
By Proposition \ref{TD}, the number $k\doteqdot\mathrm{ord\ }\hat{t}$ is
finite. For an arbitrary vector $i\in\mathbb{Z}_{+}^{n},$ $\left\vert
i\right\vert \leq m,$ we choose a test function $\varphi_{i}$ such that
$\mathrm{supp\ }\hat{\varphi}_{i}\subset B\left(  d\left(  \Sigma\right)
\right)  $ and
\[
\hat{\varphi}_{i}\left(  \xi\right)  =\frac{\left(  -1\right)  ^{\left\vert
i\right\vert }}{i!}\xi^{i}%
\]
on $B\left(  d\left(  \Sigma/2\right)  \right)  ,$ where$\ i!=i_{1}%
!...i_{n}!.$ Taking $\varphi=\varphi_{i}$ we obtain\ for any test function
$\alpha$ on $E^{\ast},$%
\[
\left[  \left(  \hat{\varphi}_{i}\ast\hat{t}\right)  \cdot\overline{\hat{t}%
}\right]  \left(  \alpha\right)  =\overline{\hat{t}}\left[  \left(
\hat{\varphi}_{i}\ast\hat{t}\right)  \left(  \alpha\right)  \right]
=\overline{\hat{t}}\left(  \xi,-D_{\xi}\right)  \left(  \hat{t}\left(
\eta,D_{\eta}\right)  \left(  \hat{\varphi}_{i}\left(  \xi-\eta\right)
\right)  \alpha\left(  \xi\right)  \right)
\]%
\[
=\left(  -1\right)  ^{\left\vert i\right\vert }\sum_{\sigma\in\Sigma}%
\overline{\hat{t}}_{\sigma}\left(  -D_{\xi}\right)  \left.  \left(  \hat
{t}_{\sigma}\left(  D_{\eta}\right)  \frac{\left(  \xi-\eta\right)  ^{i}}%
{i!}\alpha\left(  \xi\right)  \right)  \right\vert _{\eta=\xi=\sigma}.
\]
For any $\sigma\in\Sigma,$ we have%
\begin{align*}
&  \left(  -1\right)  ^{\left\vert i\right\vert }\overline{\hat{t}}_{\sigma
}\left(  -D_{\xi}\right)  \left.  \left(  \hat{t}_{\sigma}\left(  D_{\eta
}\right)  \frac{\left(  \xi-\eta\right)  ^{i}}{i!}\alpha\left(  \xi\right)
\right)  \right\vert _{\eta=\xi=\sigma}\\
&  =\left(  -1\right)  ^{\left\vert i\right\vert }\sum_{j}\overline{\hat{t}%
}_{\sigma,j}\left(  -D_{\xi}\right)  \left.  \left(  \sum_{k}\overline{\hat
{t}}_{\sigma,k}\left(  D_{\eta}\right)  \frac{\left(  \xi-\eta\right)  ^{i}%
}{i!}\alpha\left(  \xi\right)  \right)  \right\vert _{\eta=\xi=\sigma}\\
&  =\sum_{j+k=i}\overline{\hat{t}}_{\sigma j}\hat{t}_{\sigma k}\ \alpha\left(
\sigma\right)  ,
\end{align*}
where we write%
\[
\hat{t}_{\sigma}\left(  D_{\xi}\right)  =\sum_{\left\vert j\right\vert \leq
m}\hat{t}_{\sigma j}D_{\xi}^{j}.
\]
Taking the sum over $i$ with weights $\lambda^{i},$ $\lambda\in\mathbb{R}%
^{n},$ we get%
\[
\sum_{\left\vert i\right\vert \leq2m}\lambda^{i}\left[  \left(  \hat{\varphi
}_{i}\ast\hat{t}\right)  \cdot\overline{\hat{t}}\right]  \left(
\alpha\right)  =\sum_{\sigma}\sum_{i}\lambda^{i}\sum_{j+k=i}\overline{\hat{t}%
}_{\sigma j}\hat{t}_{\sigma k}\alpha\left(  \sigma\right)
\]%
\[
=\sum_{\sigma}\sum_{j}\sum_{k}\lambda^{j}\overline{\hat{t}}_{\sigma j}%
\sum\lambda^{k}\hat{t}_{\sigma k}\alpha\left(  \sigma\right)  =\sum_{\sigma
}\left\vert \sum_{k}\lambda^{k}\hat{t}_{\sigma k}\right\vert ^{2}\alpha\left(
\sigma\right)  .
\]
It follows that functional%
\[
F_{x\rightarrow\xi}\left(  A_{\varphi}\right)  =\sum_{i}F\left(  \varphi
_{i}t\ast t^{\ast}\right)  =\sum_{\sigma\in\Sigma}\left\vert \sum_{k}%
\lambda^{k}\hat{t}_{\sigma k}\right\vert ^{2}\alpha\left(  \sigma\right)
\]
is a non negative measure supported on $\Sigma$, where
\[
\varphi=\sum\lambda^{i}\varphi_{i}.
\]
We shall obtain a strictly positive measure, if choose real $\lambda$ in such
a way that the sum $\sum\lambda^{k}\hat{t}_{\sigma k}$ does not vanish for
each point $\sigma\in\Sigma.$ It can be done, since the zero set $Z_{\sigma}$
of this polynomial in $\mathbb{R}^{n}$ is nowhere dense for any $\sigma
\in\Sigma$, since $t_{\sigma}\neq0$ by the assumption.$\ $The complement of
the union $\cup Z_{\sigma}$ is also nowhere dense, since the set $\Sigma$ is
countable. Finally, for any $\lambda\in\mathbb{R}^{n}\backslash\left(  \cup
Z_{\sigma}\right)  ,\ F\left(  A_{\varphi}\right)  $ is a positive measure
supported by $\Sigma$ and the support of $A_{\varphi}$\ is contained in the
u.d. set $\Lambda-\Lambda.$ By Theorem \ref{LO} applied to $A_{\varphi},$ we
conclude that $\Sigma$ is a $\Gamma$-crystal with some lattice $\Gamma$.
$\square$

\section{Poisson structure of distributions with crystal geometry}

Now we check that for any distribution as in Theorem \ref{Gen}, the Fourier
transform can be calculated by means of the classical Poisson summation
formula%
\begin{equation}
F\left(  \sum_{p\in L}\delta_{p}\mathrm{d}x\right)  =\left\vert \det
X\right\vert \sum_{\sigma\in L^{\ast}}\delta_{\sigma}, \label{8}%
\end{equation}
where $L$ is a lattice in $E$, $L^{\ast}$ is the dual lattice and $X=\left\{
x_{i}\left(  p_{j}\right)  \right\}  $ where $p_{1},...,p_{n}$ are generators
of $L.$ Equation (\ref{8}) still holds for any shift of the lattice $L$ with
the corresponding exponential factor in the second term. For an arbitrary
linear operator $P\left(  x,D\right)  $ that is a polynomial on both its
arguments and for any vector $q\in E,$ we have
\begin{equation}
F\left(  P\left(  x,D\right)  \sum_{p\in L+q}\delta_{p}\mathrm{d}x\right)
=\exp\left(  -\mathrm{j}\left\langle \sigma,q\right\rangle \right)  P\left(
\mathrm{j}^{-1}D_{\sigma},\mathrm{j}\sigma\right)  \sum_{\sigma\in L^{\ast}%
}\delta_{\sigma}, \label{9}%
\end{equation}
where we use notation $\left\langle \sigma,q\right\rangle $ for evaluation of
the inner product of the covector $\sigma$ and the vector $q.$ Exponential
factor $\exp\left(  -\mathrm{j}\left\langle \eta,x\right\rangle \right)  $
inserted in the bracket causes the shift in the right hand side by the vector
$\eta\in E^{\ast}.$

\begin{theorem}
\label{po}Any tempered distribution $t$ supported by a\ crystal $\Lambda
=\cup_{1}^{I}\left(  L+q_{i}\right)  $ in $E$ such that $\hat{t}$ is supported
by a crystal$\ \Sigma=\cup_{1}^{K}\left(  \Gamma+\eta_{k}\right)  $ in
$E^{\ast}$ has the Poisson form%
\begin{equation}
t=\sum_{k=1}^{K}\sum_{i=1}^{I}\exp\left(  -\mathrm{j}\left\langle \eta
_{k},x\right\rangle \right)  P_{ik}\left(  x,D\right)  \sum_{p\in L+q_{i}%
}\delta_{p}, \label{16}%
\end{equation}
where $P_{ik}\left(  x,D\right)  $ are differential operators with polynomial
coefficients. Moreover, $\Gamma=L^{\ast}$ and $\hat{t}$ can be calculated by
means of (\ref{9}).
\end{theorem}

If $t$ and $\hat{t}$ are measures, all $P_{ik}$ are constants and
representation (\ref{16}) was obtained in \cite{LO15}.

\textit{Proof. }Let
\begin{equation}
G\doteqdot\left\{  x\in E;\ \left\vert \left\langle g_{i},x\right\rangle
\right\vert \leq1,\ i=1,...,n\right\}  , \label{17}%
\end{equation}
where $g_{1},...,g_{n}\in E^{\ast}$\ are generators of the lattice $L^{\ast}.$
For any $q\in L,$ we denote$\ \Lambda\left(  q\right)  =\Lambda\cap\left(
q+\left(  k+1\right)  KG\right)  ,$ where $k=\mathrm{ord\ }t.\ $Distribution%
\[
t\left[  q\right]  =\sum_{p\in\Lambda\left(  q\right)  }t_{p-q}\left(
D\right)  \delta_{p-q}%
\]
is supported by the compact set $\Lambda\left(  0\right)  $ and has order
$\leq m.$ The space $\Delta_{m}$ of all distributions supported on
$\Lambda\left(  0\right)  $ of order $\leq m$ has finite dimension. Let
$e_{1},...,e_{n}$ be the generators of $L.$ For any $i=1,...,n,$ distributions
$t\left[  je_{i}\right]  ,\ j=0,1,...,J\left(  i\right)  $ belong to
$\Delta_{m}$ and are linearly dependent$,$ if $J\left(  i\right)  $ is
sufficiently big. This yields
\[
\sum_{j=1}^{J\left(  i\right)  }c_{i}^{j}t\left[  je_{i}\right]  =0
\]
for some constants $c_{i}^{j},\ j=1,...,J\left(  i\right)  ,$ which means that
the sum%
\[
s_{i}\doteq\sum_{j}c_{i}^{j}t\left(  x-je_{i}\right)
\]
vanishes on $\Lambda\left(  0\right)  .$ The Fourier transform $\hat{\nu}_{i}$
is supported by $\Sigma,$ since so does $\hat{t}.$ This implies that $s_{i}=0$
by Proposition \ref{P} below, consequently $t$ satisfies the system of
difference equations%
\begin{equation}
\sum_{j}c_{i}^{j}t\left(  x-je_{i}\right)  =0,\ i=1,...,n. \label{13}%
\end{equation}
According to the classical result, the general solution of the system
(\ref{13}) for a scalar function $f$ on $\mathbb{Z}^{n}$ has the form%
\[
f\left(  x\right)  =\sum_{k=1}^{K}\exp\left\langle \lambda_{k}x\right\rangle
p_{k}\left(  x\right)  ,\ x\in\Lambda,
\]
where $\lambda_{k}\in\mathbb{R}^{n}$,\ $k=1,...,K$ are solutions of the
characteristic system%
\[
\sum_{j=1}c_{i}^{j}\left\langle \lambda_{k}e_{i}\right\rangle ^{j}%
=0,\ i=1,...,n,
\]
and $p_{k},$ $k=1,...,K$ are some polynomials of degree $<m_{k},$ where
$m_{k}$ is the multiplicity of $\lambda_{k}.$ The same formula holds for
solutions $t$ with values in an arbitrary vector space $V.$ For $V=\Delta_{m}%
$, we obtain%
\[
t\left(  x\right)  =\sum_{1}^{K}\exp\left\langle \lambda_{k}x\right\rangle
P_{k}\left(  x,D\right)  \delta_{x},\ x\in\Lambda,
\]
where $P_{k}$ are polynomials of $2n$ variables. Note that $\lambda_{k}$ is
pure imaginary if $P_{k}\neq0,$ since $t$ is a tempered distribution and has
at most polynomial growth at infinity. Therefore\ we can write $\lambda
_{k}=-\mathrm{j}\eta_{k}$ for a real vector $\eta_{k}.$ This implies
(\ref{9}), where $P_{ki}$ is the restriction of $P_{k}$ to the crystal
$L+q_{i},i=1,...,I:$ $\square$

\section{Gap theorem for crystals and the end of the proof}

\begin{proposition}
\label{P}Let $\hat{T}$ be a tempered distribution on $E$ such that $F\left(
T\right)  $ is supported by a $\Gamma$-crystal $\Sigma=\cup_{1}^{K}\left(
\Gamma+\eta_{k}\right)  $ in $E^{\ast}.$ If $\nu$\ vanishes on a neighborhood
of $\left(  l+1\right)  KG,$ where $G$ is as in (\ref{17}) and $\mathrm{ord}%
\ \hat{T}=l$, then $\nu=0.$
\end{proposition}

\textit{Proof. }Let $\psi$ be a test function on $E^{\ast}$ with support in
the ball $B\left(  \rho\right)  $ of small radius$\ \rho>0$ such that
$\psi\left(  0\right)  =1$. For a number $k,$ $1\leq k\leq K,$ we set%
\begin{equation}
\Psi_{k}\left(  \xi\right)  =\psi\left(  \xi\right)  \prod_{j\neq k}%
\prod_{i=1}^{n}\sin^{m+1}\left\langle h_{i},\xi-\eta_{j}\right\rangle
\prod_{i=1}^{n}\mathrm{sinc}^{m+1}\left\langle h_{i},\xi-\eta_{k}\right\rangle
, \label{6}%
\end{equation}
where $\mathrm{sinc}\lambda=\sin\lambda/\lambda$ and $h_{1},...,h_{n}$ are
generators of $\Gamma^{\ast}.$ We have $\nabla^{j}\Psi\left(  \xi\right)  =0$
for any $j=0,...,m$ and $\xi\in\Sigma\backslash\left\{  \eta_{k}\right\}  ,$
and $\nabla^{j}\Psi\left(  \eta_{k}\right)  =0,\ $for $j=1,...,m,\ \Psi\left(
\eta_{k}\right)  =1.$ It follows that $\hat{T}_{0}=\hat{T}\left(  \Psi
_{k}\right)  =T\left(  F\left(  \Psi_{k}\right)  \right)  .$ On the other
hand,
\[
\mathrm{supp\ }F\left(  \Psi_{k}\right)  \subset\left(  l+1\right)
KQ+B\left(  \rho\right)  ,
\]
if $\rho$ is sufficiently small. This follows from $\mathrm{supp\ }%
\widehat{\mathrm{sinc}}=\left[  -1,1\right]  .$ By the assumption, $T\left(
F\left(  \Psi_{k}\right)  \right)  =0$ which implies $\hat{T}_{\eta_{k},0}=0,$
since $\Psi_{k}$ vanishes on $\Sigma$ except the point $\eta_{k}=0.$ Replacing
several factors $\mathrm{sinc}\left\langle h_{i},\xi-\eta_{k}\right\rangle $
in (\ref{6}) by $\mathrm{sin}\left\langle h_{i},\xi-\eta_{k}\right\rangle ,$
we apply the same arguments and check that $\hat{T}_{\eta_{k},i}=0$ for any
$i.$ This implies that $\hat{T}_{\eta_{k}}=0.$ The same arguments can be
applied to any shift of vectors $\eta_{k}$, which gives $\hat{T}_{\sigma}%
=0$\ for any $\sigma\in\Gamma+\eta_{k}\ $and any $k.$ $\square$

\textit{End of the proof} of Theorem \ref{Gen}.\textit{ }According to Lemma
\ref{C},\ tempered function $\hat{t}=F\left(  t\right)  $ is supported by the
$\Gamma$-crystal $\Sigma,$ whereas $t$ is supported by a u.d. set $\Lambda.$
Changing roles of $E$ and $E^{\ast},$ we note set $\Sigma-\Sigma$ is u.d.,
since it is a $\Gamma$-crystal, the distribution $s\doteqdot\hat{t}%
\mathrm{d}\xi$ is supported by $\Sigma$ and distribution $t=F^{-1}\left(
s\right)  \mathrm{d}x$ is supported by u.d. set $\Lambda.$ By Lemma \ref{C},
$\Lambda$ is a $\Delta$-crystal with a lattice $\Delta$ in $E.$ By Theorem
\ref{po}, $\Gamma=\Delta^{\ast}$ which completes the proof. $\square$

\end{document}